\magnification=1200

\hsize=11.25cm    
\vsize=18cm     
\parindent=12pt   \parskip=5pt     

\hoffset=.5cm   
\voffset=.8cm   

\pretolerance=500 \tolerance=1000  \brokenpenalty=5000

\catcode`\@=11

\font\eightrm=cmr8         \font\eighti=cmmi8
\font\eightsy=cmsy8        \font\eightbf=cmbx8
\font\eighttt=cmtt8        \font\eightit=cmti8
\font\eightsl=cmsl8        \font\sixrm=cmr6
\font\sixi=cmmi6           \font\sixsy=cmsy6
\font\sixbf=cmbx6

\font\tengoth=eufm10 
\font\eightgoth=eufm8  
\font\sevengoth=eufm7      
\font\sixgoth=eufm6        \font\fivegoth=eufm5

\skewchar\eighti='177 \skewchar\sixi='177
\skewchar\eightsy='60 \skewchar\sixsy='60

\newfam\gothfam           \newfam\bboardfam

\def\tenpoint{
  \textfont0=\tenrm \scriptfont0=\sevenrm \scriptscriptfont0=\fiverm
  \def\rm{\fam\z@\tenrm}
  \textfont1=\teni  \scriptfont1=\seveni  \scriptscriptfont1=\fivei
  \def\oldstyle{\fam\@ne\teni}\let\old=\oldstyle
  \textfont2=\tensy \scriptfont2=\sevensy \scriptscriptfont2=\fivesy
  \textfont\gothfam=\tengoth \scriptfont\gothfam=\sevengoth
  \scriptscriptfont\gothfam=\fivegoth
  \def\goth{\fam\gothfam\tengoth}
  
  \textfont\itfam=\tenit
  \def\it{\fam\itfam\tenit}
  \textfont\slfam=\tensl
  \def\sl{\fam\slfam\tensl}
  \textfont\bffam=\tenbf \scriptfont\bffam=\sevenbf
  \scriptscriptfont\bffam=\fivebf
  \def\bf{\fam\bffam\tenbf}
  \textfont\ttfam=\tentt
  \def\tt{\fam\ttfam\tentt}
  \abovedisplayskip=12pt plus 3pt minus 9pt
  \belowdisplayskip=\abovedisplayskip
  \abovedisplayshortskip=0pt plus 3pt
  \belowdisplayshortskip=4pt plus 3pt 
  \smallskipamount=3pt plus 1pt minus 1pt
  \medskipamount=6pt plus 2pt minus 2pt
  \bigskipamount=12pt plus 4pt minus 4pt
  \normalbaselineskip=12pt
  \setbox\strutbox=\hbox{\vrule height8.5pt depth3.5pt width0pt}
  \let\bigf@nt=\tenrm       \let\smallf@nt=\sevenrm
  \normalbaselines\rm}

\def\eightpoint{
  \textfont0=\eightrm \scriptfont0=\sixrm \scriptscriptfont0=\fiverm
  \def\rm{\fam\z@\eightrm}
  \textfont1=\eighti  \scriptfont1=\sixi  \scriptscriptfont1=\fivei
  \def\oldstyle{\fam\@ne\eighti}\let\old=\oldstyle
  \textfont2=\eightsy \scriptfont2=\sixsy \scriptscriptfont2=\fivesy
  \textfont\gothfam=\eightgoth \scriptfont\gothfam=\sixgoth
  \scriptscriptfont\gothfam=\fivegoth
  \def\goth{\fam\gothfam\eightgoth}
  
  \textfont\itfam=\eightit
  \def\it{\fam\itfam\eightit}
  \textfont\slfam=\eightsl
  \def\sl{\fam\slfam\eightsl}
  \textfont\bffam=\eightbf \scriptfont\bffam=\sixbf
  \scriptscriptfont\bffam=\fivebf
  \def\bf{\fam\bffam\eightbf}
  \textfont\ttfam=\eighttt
  \def\tt{\fam\ttfam\eighttt}
  \abovedisplayskip=9pt plus 3pt minus 9pt
  \belowdisplayskip=\abovedisplayskip
  \abovedisplayshortskip=0pt plus 3pt
  \belowdisplayshortskip=3pt plus 3pt 
  \smallskipamount=2pt plus 1pt minus 1pt
  \medskipamount=4pt plus 2pt minus 1pt
  \bigskipamount=9pt plus 3pt minus 3pt
  \normalbaselineskip=9pt
  \setbox\strutbox=\hbox{\vrule height7pt depth2pt width0pt}
  \let\bigf@nt=\eightrm     \let\smallf@nt=\sixrm
  \normalbaselines\rm}

\tenpoint

\def\pc#1{\bigf@nt#1\smallf@nt}         \def\pd#1 {{\pc#1} }

\frenchspacing

\def\raggedbottom{\topskip 10pt plus 36pt\r@ggedbottomtrue}

\def\pointir{\unskip . --- \ignorespaces}

\def\Medbreak{\vskip-\lastskip\medbreak}

\long\def\th#1 #2\enonce#3\endth{
   \Medbreak\noindent
   {\pc#1} {#2\unskip}\pointir{\it #3}\smallskip}

\def\decale#1{\smallbreak\hskip 28pt\llap{#1}\kern 5pt}
\def\decaledecale#1{\smallbreak\hskip 34pt\llap{#1}\kern 5pt}
\def\puce{\smallbreak\hskip 6pt{$\scriptstyle\bullet$}\kern 5pt}

\def\eqalign#1{\null\,\vcenter{\openup\jot\m@th\ialign{
\strut\hfil$\displaystyle{##}$&$\displaystyle{{}##}$\hfil
&&\quad\strut\hfil$\displaystyle{##}$&$\displaystyle{{}##}$\hfil
\crcr#1\crcr}}\,}

\catcode`\@=12

\showboxbreadth=-1  \showboxdepth=-1

\newcount\numerodesection \numerodesection=1
\def\section#1{\bigbreak
 {\bf\number\numerodesection.\ \ #1}\nobreak\medskip
 \advance\numerodesection by1}

\mathcode`A="7041 \mathcode`B="7042 \mathcode`C="7043 \mathcode`D="7044
\mathcode`E="7045 \mathcode`F="7046 \mathcode`G="7047 \mathcode`H="7048
\mathcode`I="7049 \mathcode`J="704A \mathcode`K="704B \mathcode`L="704C
\mathcode`M="704D \mathcode`N="704E \mathcode`O="704F \mathcode`P="7050
\mathcode`Q="7051 \mathcode`R="7052 \mathcode`S="7053 \mathcode`T="7054
\mathcode`U="7055 \mathcode`V="7056 \mathcode`W="7057 \mathcode`X="7058
\mathcode`Y="7059 \mathcode`Z="705A


\def\hfl#1#2#3{\smash{\mathop{\hbox to#3{\rightarrowfill}}\limits
^{\textstyle#1}_{\textstyle#2}}}

\def\ogoth{{\goth o}}

\def\pgoth{{\goth p}}

\def\Q{{\bf Q}}

\def\N{{\bf N}}

\def\Z{{\bf Z}}

\def\F{{\bf F}}

\def\Hom{\mathop{\rm Hom}\nolimits}

\def\Card{\mathop{\rm Card}\nolimits}
\def\Gal{\mathop{\rm Gal}\nolimits}

\def\series#1{(\!(#1)\!)}

\def\to{\rightarrow}
\def\lcm{\mathop{\rm lcm}\nolimits}

\def\mod{\mathop{\rm mod.}\nolimits}
\def\pmod#1{\;(\mod#1)}

\def\Aut{\mathop{\rm Aut}\nolimits}

\def\boxit#1{\vbox{\hrule\hbox{\vrule\kern1pt
       \vbox{\kern1pt#1\kern1pt}\kern1pt\vrule}\hrule}}
\def\cqfd{\hfill\boxit{\phantom{\i}}}

\newcount\numero \numero=1
\def\numeroter{{({\oldstyle\number\numero})}\ \advance\numero by1}

\newcount\refno 
\long\def\ref#1:#2<#3>{                                        
\global\advance\refno by1\par\noindent                              
\llap{[{\bf\number\refno}]\ }{#1} \pointir{\it #2} #3\goodbreak }

\newcount\refno 
\long\def\ref#1:#2<#3>{                                        
\global\advance\refno by1\par\noindent                              
\llap{[{\bf\number\refno}]\ }{#1} \pointir{\it #2} #3\goodbreak }

\def\citer#1(#2){[{\bf\number#1}\if#2\empty\relax\else,\ {#2}\fi]}

\def\boxit#1{\vbox{\hrule\hbox{\vrule\kern1pt
       \vbox{\kern1pt#1\kern1pt}\kern1pt\vrule}\hrule}}
\def\cqfd{\hfill\boxit{\phantom{\i}}}

\newbox\bibbox
\setbox\bibbox\vbox{\bigbreak
\centerline{{\pc BIBLIOGRAPHY}}

\ref{\pc BOURBAKI} (N):
Alg\`ebre, Chapitres 4 \`a 7,
<Masson, Paris, 1981, 422~pp.>
\newcount\bourbaki \global\bourbaki=\refno

\ref{\pc DALAWAT} (C):
Further remarks on local discriminants,
<J.\ Ramanujan Math.\ Soc. {bf 25} (2010) 4, 393--417.
Cf.~arXiv\string:0909.2541.>
\newcount\further \global\further=\refno

\ref{\pc DALAWAT} (C):
Serre's ``\thinspace formule de masse\thinspace'' in prime degree,
<Monats\-hefte Math.\ {\bf 166} (2012) 1, 73--92.
Cf.~arXiv\string:1004.2016.>
\newcount\monatshefte \global\monatshefte=\refno

\ref{\pc DEL \pc CORSO} (I), {\pc DVORNICICH} (R) \& {\pc MONGE} (M):
On wild extensions of a $p$-adic field,
<J.\ Number Theory {\bf 174} (2017), 322--342.  Cf.~aXiv\string:1601.05939.> 
\newcount\deldvomonge \global\deldvomonge=\refno

\ref{\pc IWASAWA} (K):
On Galois groups of local fields.
<Trans.\ Amer.\ Math.\ Soc.\ {\bf 80} (1955), 448--469.>
\newcount\iwasawa \global\iwasawa=\refno

\ref{\pc JANNSEN} (U):
\"Uber Galoisgruppen lokaler K\"orper.
<Invent.\ Math.\ {\bf 70} (1982/83) 1, 53--69.>
\newcount\jannsen \global\jannsen=\refno

\ref{\pc KOCH} (H):
Galois theory of p-extensions. 
<Springer-Verlag, Berlin, 2002, 190~pp.>
\newcount\kochpext \global\kochpext=\refno




} 

\centerline{\bf Little galoisian modules} 
\bigskip\bigskip 
\centerline{Chandan Singh Dalawat} 
\centerline{Harish-Chandra Research Institute}
\centerline{Chhatnag Road, Jhunsi, Allahabad 211019, India} 
\centerline{\tt dalawat@gmail.com}

\bigskip\bigskip

{{\bf Abstract}.  Let $p$ be a prime number, let $K$ be a $p$-field (a
  local field with finite residue field of characteristic~$p$), let
  $L$ be a finite galoisian tamely ramified extension of $K$, and let
  $G=\Gal(L|K)$. Suppose that $L$ is split over $K$ in the sense that
  the short exact sequence $1\to T\to G\to G/T\to1$ has a section,
  where $T$ is the inertia subgroup of $G$.  We determine the
  structure of the $\F_p[G]$-module $L^\times\!/L^{\times p}$ in
  characteristic~$0$ when the $p$-torsion subgroup ${}_pL^\times$ of
  $L^\times$ has order~$p$, and of the $\F_p[G]$-modules
  $L^\times\!/L^{\times p}$ and $L^+\!/\wp(L^+)$ in
  characteristic~$p$, where $\wp(x)=x^p-x$.

Let $\tilde K$ be a maximal galoisian extension of $K$, let $V$ be the
  maximal tamely ramified extension of $K$ in $\tilde K$, let
  $\Gamma=\Gal(V|K)$, and let $B$ be the maximal abelian extension of
  exponent~$p$ of $V$ in $\tilde K$.  We determine the structure of
  the $\F_p[[\Gamma]]$-module $\Gal(B|V)$, and show how this leads in
  characteristic~$0$ to a simple proof of the fact that the profinite
  group $\Gal(\tilde K|K)$ is generated by $[K:\Q_p]+3$ elements.

\footnote{}{{\it MSC2010~:} Primary 11R23, 11S15}
\footnote{}{{\it Keywords~:} Local fields, galoisian modules, tame ramification 
}}

\bigskip\bigskip

{\bf 1.  Introduction}
\bigskip

\numeroter   In the first part of this Note~(\S2),
we work with a finite field $k$, a finite extension $l$ of $k$, and an
injective morphism of groups $\theta:T\to l^\times$.  Let $q=\Card k$,
let $\Sigma=\Gal(l|k)$, and let $\sigma$ be the generator $x\mapsto
x^q$ ($x\in l$) of $\Sigma$.  View $T$ as a submodule of the
$\Sigma$-module $l^\times$, and let $G=T\times_q\Sigma$ be the twisted
product of $\Sigma$ by $T$.  For every $i\in\Z$, we have the
$k[G]$-module $l(i)$ whose underlying $k$-space is $l$ and on which
$G$ acts by
$$
\sigma.x=x^q,\quad t.x=\theta(t)^ix\qquad (x\in l,\ t\in T),
$$
We show that these modules are projective, and determine when two such
modules are isomorphic.  The main tool is a lemma of
Iwasawa \citer\iwasawa(Lemma~1, p.~449).

\numeroter In the second and third parts (\S3 and \S4), we work with a
local field $K$ with finite residue field $k$ of cardinality $q$ and
characteristic~$p$, and a finite galoisian tamely ramified split
extension $L$ of $K$ of residue field $l$ and group $G=\Gal(L|K)$.
The inertia subgroup $T\subset G$ comes with a faithful character
$\theta:T\to l^\times$ and, since $L$ is split over $K$ by hypothesis,
$G$ is isomorphic to $T\times_q\Sigma$, where $\Sigma=\Gal(l|k)$.

Exploiting our study of the $k[G]$-modules $l(i)$ in \S2, we determine
the structure of the $\F_p[G]$-module $L^\times\!/L^{\times p}$ when
$K$ has characteristic~$0$ and ${}_pL^\times$ has order~$p$ in \S3,
and of the $\F_p[G]$-modules $L^\times\!/L^{\times p}$ and
$L^+\!/\wp(L^+)$, where $\wp(x)=x^p-x$ ($x\in L$), when $K$ has
characteristic~$p$ in \S4.  The results are a generalisation from the
much simpler case $L=K(\!\root p-1\of{K^\times})$ treated
in \citer\monatshefte().

\numeroter In the fourth part (\S5 and \S6), we consider the maximal
tamely ramified extension $V$ of $K$ of group $\Gamma=\Gal(V|K)$ and,
putting everything together, determine the structure of the
$\F_p[[\Gamma]]$-module $V^\times\!/V^{\times p}$ in
characteristic~$0$ and of the $\F_p[[\Gamma]]$-modules
$V^\times\!/V^{\times p}$ and $V^+\!/\wp(V^+)$ in characteristic~$p$.
As a consequence, we determine (in both cases~: mixed- and
equi-characteristic) the structure of the $\F_p[[\Gamma]]$-module
$\Gal(B|V)$, where $B$ is the maximal abelian extension of $V$ of
exponent~$p$.  This is achieved by passing to the limit over the
results of \S3 and \S4.

Finally, in \S6 we take $K$ to be a finite extension of $\Q_p$ with
maximal galoisian extension $\tilde K$ and show how the sturcture
theorem for the $\F_p[[\Gamma]]$-module $\Gal(B|V)$ as proved in \S5
leads to a simple proof of the fact that the profinite group
$\Gal(\tilde K|K)$ is generated by $[K:\Q_p]+3$ elements.

\bigbreak
{\bf 2. Iwasawa's lemma}
\bigskip

\numeroter  We recall a crucial lemma from Iwasawa \citer\iwasawa() and
simplify its proof.  Let $p$ be a prime number and let $e>0$ be an
integer $\not\equiv0\pmod p$.  Let $g>0$ be a multiple of the order of
$\bar p\in(\Z/e\Z)^\times$, so that there is a unique morphism of
groups $\Z/g\Z\to(\Z/e\Z)^\times$ such that $1\mapsto p$.  Let $n$ be
a multiple of $\lcm(p-1,e)$, and write $n=c.(p-1)$ and $n=d.e$.  Let
$b^{(i)}$ ($i>0$) be the sequence of positive integers
$\not\equiv0\pmod p$, namely $b^{(i)}=i+\lfloor(i-1)/(p-1)\rfloor$.
For every $a,b\in\Z$, we denote by $[a,b]$ the set of {\it integers\/}
between $a$ and $b$.

\numeroter {\it The map\/ $[1,n]\times\Z/g\Z\to\Z/e\Z$ sending\/
  $(i,j)$ to\/ $b^{(i)}p^j\pmod e$ is surjective and every fibre has\/
  $dg$ elements.}

{\it Proof}.  Consider the map $[1,cp]\times\Z/g\Z\to\Z/e\Z$ sending
$(r,j)$ to $rp^j\pmod e$~; it is the ``\thinspace product\thinspace''
of the natural map of reduction $\pmod e$ on the first factor and the
map $1\mapsto p$ discussed in ($\oldstyle4$) on the
second factor, and it is clearly surjective.  View the interval
$[1,cp]$ as the (disjoint) union of the successive intervals $[1,de]$
and $[n+1,n+c]$ on the one hand, and as the (disjoint) union of the
subsets $(b^{(i)})_{i\in[1,n]}$ and $(ip)_{i\in[1,c]}$ on the other.
By the {\it contribution\/} of a subset $S\subset[1,cp]$ we mean the
family $(t_x)_{x\in\Z/e\Z}$, where $t_x$ is the number of antecedents
of $x$ in $S\times\Z/g\Z$.  Clearly, the contribution of
$(ip)_{i\in[1,c]}$ is the same as that of $[n+1,n+c]$, because
$j\mapsto j+1$ is a permutation of $\Z/g\Z$.  So the contribution of
$(b^{(i)})_{i\in[1,n]}$ is the same as that of $[1,de]$, which is easy
to compute~: for fixed $x\in\Z/e\Z$ and $j\in\Z/g\Z$, there are
exactly $d$ elements $r\in[1,de]$ such that $rp^j\equiv x\pmod
e$.  \cqfd

\numeroter {\it The\/ $k[\Sigma]$-module\/ $l$}.  Let $k$ be a finite
extension of $\F_p$ of cardinality $q=p^a$.  Let $l$ be a finite
extension of $k$, and put $f=[l:k]$, $g=af$.  Let $\Sigma=\Gal(l|k)$,
and denote by $\sigma$ the generator $x\mapsto x^q$ ($x\in l$) of
$\Sigma$.  The $k[\Sigma]$-module $l$ is free of rank~$1$, as follows
from the normal basis theorem
\citer\bourbaki(V.70)~: there exists an $\alpha\in l$ such that
$(\sigma^i(\alpha))_{i\in\Z/f\Z}$ is a $k$-basis of $l$.

\numeroter {\it The groups\/ $T$, $G$, and the character\/
  $\theta:T\to l^\times$}. Let $T$ be a subgroup of $l^\times$, $e$
its order (so that $q^f\equiv1\pmod e$) and $\theta:T\to l^\times$ the
inclusion~; the group $\Hom(T,l^\times)$ of characters of $T$ is
cyclic of order~$e$ and generated by~$\theta$. Identifying $\Aut(T)$
with $(\Z/e\Z)^\times$, there is a unique morphism of groups
$\Sigma\to\Aut(T)$ such that $\sigma\mapsto q$~; endow $T$ with this
action of $\Sigma$ (which is the galoisian action as a subgroup of
$l^\times$) and let $G=T\times_q\Sigma$ be the twisted product of
$\Sigma$ by the $\Sigma$-module $T$, so that $\sigma t\sigma^{-1}=t^q$
for every $t\in T$~; we sometimes write $G=T\Sigma$.

\numeroter Concretely, if we choose a generator $\tau$ for $T$, then
the group $G$ has the presentation $G=\langle\sigma,\tau\mid
\sigma^f=1,\tau^e=1,\sigma\tau\sigma^{-1}=\tau^q\rangle$, and
$\theta(\tau)$ is a primitive $e$-th root of~$1$ in $l$.  Conversely,
if we choose an element $\eta\in l^\times$ of order~$e$, then
$\theta^{-1}(\eta)$ is a generator of $T$.  In what follows, we don't
need to choose $\tau$ or $\eta$.

\numeroter {\it The\/ $k[G]$-modules\/ $l(r)$}. For every $r\in\Z$,
make $G$ act on the $k$-space $l$ by the law
$$
\sigma.x=x^q,\quad t.x=\theta(t)^rx,\qquad (x\in l,\ t\in T)
$$
and denote the resulting $k[G]$-module by $l(r)$~; it is clear that
$l(r)$ depends only on the image $\bar r\in\Z/e\Z$, and that the
$k[G]$-module $l(0)$ is deduced from the $k[\Sigma]$-module $l$ via
the map $k[G]\to k[\Sigma]$ coming from the projection $G\to\Sigma$.

\numeroter  {\it The\/ $k[G]$-module $\bigoplus\limits_{i\in\Z/e\Z} l(i)$ is
free of rank\/~$1$.}

{\it Proof}.  Indeed, the $k[\Sigma]$-module $l$ is free of rank~$1$
$\oldstyle(6)$, and if $\alpha\in l$ is a $k[\Sigma]$-basis of $l$,
then $\alpha$ is a $k[G]$-basis of $\bigoplus\limits_{i\in\Z/e\Z}
l(i)$, where the $k[G]$-module $l(0)$ has been identified with
$l$ as in $\oldstyle(9)$.  \cqfd

\numeroter {\it The\/ $l[G]$-module\/ 
$l(r)\bigotimes\limits_{\F_p[G]}l[G]$ is isomorphic to\/
$\bigoplus\limits_{j\in\Z/g\Z}\theta^{rp^j}$ on which $\sigma$ acts by\/
$(x_j)_{j\in\Z/g\Z}\mapsto(x_{j+a})_{j\in\Z/g\Z}$, where\/ $a=[k:\F_p]$
and\/ $g=af$}.

{\it Proof}.  Here we have abused notation to make
$\chi\in\Hom(T,l^\times)$ stand for an $l$-line on which $T$ acts via
$\chi$.  By the normal basis theorem \citer\bourbaki(V.70), there
exists a $\beta\in l$ such that the $\beta_j=\beta^{p^j}$
($j\in\Z/g\Z$) constitute an $\F_p$-basis of $l$.  When we fix such a
$\beta$, we get a $l$-basis $\gamma_j=\beta_j\otimes1$ of
$l(r)\otimes_{\F_p[G]}l[G]$.  The $l$-linear actions of $\sigma$ and $T$ on
the $\gamma_j$ are given by
$$
\sigma.\gamma_j=(\beta^{p^j})^q\otimes1=\beta^{p^{j+a}}\otimes1=\gamma_{j+a}
$$
and
$$
t.\gamma_j=t.(\beta^{p^j}\otimes1)
=(t.\beta)^{p^j}\otimes1
=\theta(t)^{rp^j}\gamma_j\qquad
(t\in T).
$$
In other words, $l(r)\otimes_{\F_p[G]}l[G]$ is isomorphic to the direct sum
of the characters $\theta^{rp^j}$ ($j\in\Z/g\Z$) which are
permuted by $\sigma$ according to $j\mapsto j+a$. \cqfd

\numeroter The group $\Gal(l|\F_p)$ acts on the set $\Hom(T,l^\times)$
of characters of $T$~: the generator $\varphi:x\mapsto x^p$ ($x\in l$)
of $\Gal(l|\F_p)$ sends a character $\chi\in\Hom(T,l^\times)$ to
$\chi^p$.

\numeroter {\it The\/ $\F_p[G]$-modules\/ $l(r)$ and\/ $l(s)$ are
  isomorphic if and only if the characters\/
  $\theta^r,\theta^s\in\Hom(T,l^\times)$ are in the same
  $\varphi$-orbit.}

{\it Proof}.  The $\F_p[G]$-modules $l(r)$ and $l(s)$ are isomorphic
if and only if the $l[G]$-modules $l(r)\otimes_{\F_p[G]}l[G]$ and
$l(s)\otimes_{\F_p[G]}l[G]$ are isomorphic \citer\bourbaki(V.70).  The
result then follows from the explicit description ($\oldstyle11$) of the latter
modules.  \cqfd

\numeroter {\it Let\/ $n$ be a multiple of\/ $\lcm(p-1,e)$,
  and write $n=de$.  The\/ $\F_p[G]$-module\/
  $M=\bigoplus\limits_{i\in[1,n]}l(b^{(i)})$
  is isomorphic to\/ $k[G]^d$.}

{\it Proof}.  It is enough to show that the $l[G]$-modules
$M\otimes_{\F_p[G]}l[G]$ and $k[G]^d\otimes_{\F_p[G]}l[G]$ are
isomorphic \citer\bourbaki(V.70).  This follows from the description
($\oldstyle11$) of these modules, the criterion ($\oldstyle13$) for
$l(r)$ and $l(s)$ to be $\F_p[G]$-isomorphic, and the numerical lemma
($\oldstyle5$).  \cqfd

\numeroter {\it The\/ $k[G]$-module\/ (resp.\/~$\F_p[G]$-module) $l(r)$
is projective.}

{\it Proof}.  This follows from the fact that $l(r)$ is a
direct summand ($\oldstyle10$) of the free module $k[G]$ (of rank~$1$ over
$k[G]$ and rank $a$ over $\F_p[G]$).  \cqfd

\numeroter (Iwasawa, \citer\iwasawa(p.~449)) {\it Let\/ $n$ be a
   multiple of\/ $\lcm(p-1,e)$, and write\/ $n=de$.  Every\/
  $\F_p[G]$-module\/ $M$ endowed with a filtration
$$
\{0\}=M_{n+1}\subset M_n\subset\cdots\subset M_2\subset M_1=M
$$
such that\/ $M_i/M_{i+1}$ is isomorphic to\/ $l(b^{(i)})$ for\/
$i\in[1,n]$ is isomorphic to\/ $k[G]^d$.}

{\it Proof}.  As the $\F_p[G]$-modules $l(r)$ are projective ($\oldstyle15$),
the filtration on $M$ splits in the sense that $M$ is
$\F_p[G]$-isomorphic to $\bigoplus\limits_{i\in[1,n]}l(b^{(i)})$.  But this
module is isomorphic to $k[G]^d$, as we have seen in ($\oldstyle14$).  \cqfd

\bigbreak
{\bf 3.  The mixed-characteristic case}
\bigskip

\numeroter Let $K$ be a finite extension of $\Q_p$ of ramification
index $e_K$ and residual degree $f_K$.  Let $L$ be a finite tamely
ramified galoisian extension of $K$ of ramification index $e$ and
residual degree $f$.  Let $L_0$ be the maximal unramified extension of
$K$ in $L$ and let $\Sigma=\Gal(L_0|K)$, $T=\Gal(L|K_0)$, and
$G=\Gal(L|K)$.  Suppose that $L$ is {\it split\/} over $K$ in the
sense that the short exact sequence $1\to T\to G\to\Sigma$ has a
section.  Equivalently \citer\iwasawa(2.1), $L=L_0(\root e\of{\pi_K})$
for some uniformiser $\pi_K$ of $K$.

Let $k$ (resp.~$l$) be the residue field of $K$ (resp.~$L$).  Identify
$\Sigma$ with $\Gal(l|k)$, and let $\sigma$ be the generator $x\mapsto
x^{q}$ ($x\in l$, $q=p^{f_K}$) of $\Sigma$, so that
$\sigma=\varphi^{f_K}$, in the notation of $\oldstyle(12)$.  The
inertia group $T$ comes equipped with a (faithful) character
$\theta:T\to l^\times$ giving the action of $T$ on the set of $e$-th
roots of $\pi_K$, where $\pi_K$ is any uniformier of $K$ such that
$L=L_0(\root e\of{\pi_K})$.  Thus, we are in the situation described
in ($\oldstyle7$).

\numeroter Assume that the $p$-torsion subgroup ${}_pL^\times$ of 
$L^\times$ has order~$p$, so that the abslolute ramification index $e_L$
of $L$ is divisible by $p-1$~; write $e_L=c.(p-1)$.  Note that every
finite tamely ramified extension of $K$ is contained in a finite
galoisian tamely ramified split extension of $K$ containing a
primitive $p$-th root of~$1$.

\numeroter For $r>0$, denote by $\bar U_L^r$ the $G$-submodule the image
of $U_L^r$ in $\overline{L^\times}=L^\times\!/L^{\times p}$, where
$U_L^r=1+\pgoth_L^r$ and $\pgoth_L$ is the unique maximal ideal of the
ring of integers $\ogoth_L$ of $L$.  Note that $\bar U_L^1$ is equal
to the image $\overline{\ogoth_K^\times}$ of the group
$\ogoth_L^\times$ of units of $\ogoth_L$.  It is known that $\bar
U_L^{r}=\{1\}$ for $r>cp$, that $\bar U_L^{cp}$ is
$\F_p[G]$-isomorphic to ${}_pL^\times$ (but the submodule $\bar
U_L^{cp}\subset\overline{L^\times}$ is not to be confused with the
submodule ${}_pL^\times\subset L^\times$ nor with its image
$\overline{{}_pL^\times}\subset\overline{L^\times}$), that $\bar
U_L^{ip}=\bar U_L^{ip+1}$ for $i\in[1,c[$ and that $\bar U_L^r/\bar
U_L^{r+1}$ is isomorphic to $l(r)$ for $r\in[1,cp[$,
$r\not\equiv0\pmod p$, where $l(r)$ is the $k[G]$-module $l$ with
$\sigma$ acting by $x\mapsto x^{q}$ and $T$ acting by the character
$\theta^r$, as in ($\oldstyle9$).  See for
example \citer\monatshefte(), where we had
$L=K(\!\root{p-1}\of{K^\times})$ but the same proofs work without
change.

\numeroter Take $n=e_L$, which is a multiple of $\lcm(p-1,e)$
($\oldstyle18$), as required in ($\oldstyle16$), and note that
$cp=b^{(n)}+1$.  For every $i\in[1,n]$, put $M_i=\bar
U_L^{b^{(i)}}/\bar U_L^{cp}$ and put $M_{n+1}=\{1\}$.  This filtration
on the $\F_p[G]$-module $M=M_1$ has the property that $M_i/M_{i+1}$ is
isomorphic to $l(b^{(i)})$ for every $i\in[1,n]$, as we have recalled
($\oldstyle19$).  It follows from Iwasawa's lemma ($\oldstyle16$) that
$M$ is isomorphic to $k[G]^{e_K}$, and hence $\bar U_L^1$ is
isomorphic to ${}_pL^\times\oplus k[G]^{e_K}$.  In summary, we get the
following result of Iwasawa \citer\iwasawa(p.~461).

\numeroter {\it Let\/ $K$ be a finite extension of\/ $\Q_p$ of
ramification index\/~$e_K$ and residue field\/ $k$, and let\/ $L$ be a
finite galoisian tamely ramified split extension of\/ $K$ of group\/
$G=\Gal(L|K)$ such that\/ ${}_pL^\times$ has order\/~$p$.  The\/
$\F_p[G]$-module\/ $\bar U_L^1$ is isomorphic to\/ ${}_pL^\times\oplus
k[G]^{e_K}$, which is isomorphic to\/
${}_pL^\times\oplus\F_p[G]^{[K:\Q_p]}$.} \cqfd

\numeroter Let $\bar\F_p$ be the maximal galoisian extension of
$\F_p$.  As a corollary, we deduce that the $\bar\F_p[G]$-module
$M\otimes_{\F_p[G]}\bar\F_p[G]$ is isomorphic to
$\bar\F_p[G]^{[K:\Q_p]}$, thereby
recovering \citer\deldvomonge(Corollary 4.7).  Our method also shows
that if $L$ is a finite galoisian tamely ramified split extension of
$K$ such that ${}_pL^\times$ is trivial but $e_L$ is divisible by
$p-1$, then the $\F_p[G]$-module $\bar
U_L^1=\overline{\ogoth_L^\times}$ is isomorphic to $k[G]^{e_K}$ and to
$\F_p[G]^{[K:\Q_p]}$, for the only difference in this case is that
$\bar U_L^{cp}$ is trivial.

\numeroter  As a curiosity, the reader may wish to determine the
structure of the $\F_p[G]$-modules
$\prod\limits_{i>0}U_L^{b^{(i)}}\!/U_L^{b^{(i)}+1}$ and
$\bigoplus\limits_{i>0}\pgoth_L^{-b^{(i)}}\!/\pgoth_L^{-b^{(i)}+1}$,
where $b^{(i)}$ is the sequence of positive integers $\not\equiv0\pmod
p$, as throughout.

\bigbreak
{\bf 4.  The equi-characteristic case}
\bigskip

\numeroter Let $K$ be a local field of characteristic~$p$ with finite
resiude field $k$ of cardinality $q$, and let $L$ be a finite
galoisian tamely ramified split extension of $K$ of ramification index
$e$ ($\not\equiv0\pmod p$) and residual degree $f$ (so that
$q^f\equiv1\pmod e$).  Every finite tamely ramified extension of $K$
is contained in such an $L$.  Concretely, the residue field $l$ of $L$
is the finite extension of $k$ of degree~$f$ and there is a
uniformiser $\pi_K$ of $K$ such that $K=k\series{\pi_K}$ and
$L=l\series{\root e\of{\pi_K}}$.  The groups $\Sigma=\Gal(l|k)$,
$G=\Gal(L|K)$, $T=\Gal(L|l\series{\pi_K})$ and the character
$\theta:T\to l^\times$ giving the action of $T$ on the set of $e$-th
roots of $\pi_K$ have the properties required in $\oldstyle(7)$.

\numeroter Let us determine the structure
of the $\F_p[G]$-module $L^+\!/\wp(L^+)$, where $\wp(x)=x^p-x$ ($x\in
L$).  It will turn out that $L^+\!/\wp(L^+)$ is isomorphic to
$\F_p\oplus k[G]^{(\N)}$, just as in the case $e=p-1$, $f=p-1$ treated
in \citer\monatshefte().  The proof combines ideas
from \citer\monatshefte() with Iwasawa's lemma ($\oldstyle16$), and is analogous
to the proof of ($\oldstyle21$).

\numeroter Let $\pgoth_L$ be the unique maximal ideal of the ring of
integers $\ogoth_L=l[[\root e\of{\pi_K}]]$ of $L$.  For every
$r\in\Z$, denote by $\overline{\pgoth_L^r}$ the $G$-submodule the
image of $\pgoth_L^r$ in $\overline{L^+}=L^+\!/\wp(L^+)$.  It is known
that $\overline{\pgoth_L^r}=\{0\}$ for $r>0$, that
$\overline{\ogoth_L^+}=\overline{\pgoth_L^0}$ is canonically
isomorphic to $\F_p$, that
$\overline{\pgoth_L^{ip+1}}=\overline{\pgoth_L^{ip}}$ for all $i<0$,
and that $\overline{\pgoth_L^r}/\overline{\pgoth_L^{r+1}}$ is
isomorphic to $l(r)$ for every $r<0$, $r\not\equiv0\pmod p$, in the
notation of ($\oldstyle9$).  See for example \citer\further() for the
general case and \citer\monatshefte() for the special case
$L=K(\!\root p-1\of{K^\times})$.

\numeroter Let $n$ be a multiple of $\lcm(p-1,e)$
($\oldstyle16$), write $n=c.(p-1)$, $n=d.e$ and note that
$cp=b^{(n)}+1$.  For $i\in[1,n]$, put $\lambda(i)=b^{(i)}-cp$ and
define $M_i=\overline{\pgoth_L^{\lambda(i)}}/\overline{\pgoth_L^0}$~;
also put $M_{n+1}=\{0\}$.  We thus get a filtration on the
$\F_p[G]$-module $M=M_1$ such that $M_i/M_{i+1}$ is isomorphic to
$l(b^{(i)}-c)$ for every $i\in[1,n]$, since $\lambda(i)\equiv
b^{(i)}-c\pmod e$.  Replacing $n$ by a suitable multiple of $n$, we
may assume that $c\equiv0\pmod e$ and therefore that $M_i/M_{i+1}$ is
isomorphic to $l(b^{(i)})$. Since $M$ is isomorphic to $k[G]^d$
($\oldstyle16$), $\overline{\pgoth_L^{-b^{(n)}}}$ is
$\F_p[G]$-isomorphic to $\F_p\oplus k[G]^d$ ($\oldstyle26$).

\numeroter Replacing $n$ by $mn$ ($m>0$), we conclude that
$\overline{\pgoth_L^{-b^{(mn)}}}$ is $\F_p[G]$-isomorphic to
$\F_p\oplus k[G]^{md}$.  As $L^+$ is the direct limit of the
$\pgoth_L^{-b^{(mn)}}$ when $m\to+\infty$, we conclude that
$L^+\!/\wp(L^+)$ is isomorphic to $\F_p\oplus k[G]^{(\N)}$, just as in
the case $L=K(\!\root p-1\of{K^\times})$ treated
earlier \citer\monatshefte().  Let us summarise.

\numeroter {\it Let\/ $K$ be local field of characteristic\/ $p$ with
finite residue field\/ $k$, let\/ $L$ be a finite galoisian tamely
ramified split extension of\/ $K$, and let\/ $G=\Gal(L|K)$.  The\/
$\F_p[G]$-module\/ $\overline{L^+}=L^+\!/\wp(L^+)$ is isomorphic to\/
$\F_p\oplus k[G]^{(\N)}$ and to\/ $\F_p\oplus\F_p[G]^{(\N)}$.} \cqfd
 
\numeroter While we are at it, we might as well determine the
structure of the $\F_p[G]$-module $\bar U_L^1$, where $\bar U_L^r$ ($r>0$)
is the image of $U_L^r=1+\pgoth_L^r$ in $L^\times\!/L^{\times p}$.  It
is easy to see that $\bar U_L^{ip}=\bar U_L^{ip+1}$ for every $i>0$ and
that $\bar U_L^r/\bar U_L^{r+1}$ is $\F_p[G]$-isomorphic to $l(r)$ for
every $r>0$, $r\not\equiv0\pmod p$.

\numeroter Let $n$ be a multiple of $\lcm(p-1,e)$ ($\oldstyle14$), with
$n=c.(p-1)$, $n=d.e$.  For $i\in[1,n]$, put $M_i=\bar U_L^{b^{(i)}}/\bar
U_L^{cp}$, and put $M_{n+1}=\{1\}$.  This filtration on the module
$M=M_1$ has the properties required for applying Iwasawa's lemma ($\oldstyle16$),
therefore $\bar U_L^1/\bar U_L^{cp}$ is isomorphic to $k[G]^d$.  Replacing
$n$ by a multiple $mn$, we see that $\bar U_L^1/\bar U_L^{mcp}$ is
isomorphic to $k[G]^{md}$.  Taking the projective limit as
$m\to+\infty$, we get the structure of $\bar U_L^1$ :

\numeroter {\it The image\/ $\bar U_L^1$ of\/ $U_L^1=1+\pgoth_L$ in\/
$L^\times\!/L^{\times p}$ is\/ $\F_p[G]$-isomorphic to\/ $k[G]^{\N}$
and to\/ $\F_p[G]^{\N}$.} \cqfd

\bigbreak
{\bf 5. Passing to the tame limit}
\bigskip

\numeroter Let $K$ be a local field with finite residue field $k$ of
characteristic~$p$ and cardinality $q$ , $V$ the maximal tamely
ramified extension of $K$, and $B$ the maximal abelian extension of
$V$ of exponent $p$, so that $B=V(\!\root p\of{V^\times})$ if $K$ has
characteristic~$0$ and $B=V(\wp^{-1}(V))$ if $K$ has
characteristic~$p$.  The pro-$p$-group $\Gal(B|V)$ is an
$\F_p[[\Gamma]]$-module, where $\Gamma=\Gal(V|K)$, and we would like
to determine its structure.  This is achieved by studying the dual
$\F_p[[\Gamma]]$-module, namely $V^\times\!/V^{\times p}$ in
characteristic~$0$ and $V^+\!/\wp(V^+)$ in characteristic~$p$.

\numeroter If $K$ has characteristic~$0$, there is an interesting
intermediate extension $B'$ which may be called the maximal {\it peu
ramifi\'ee\/} extension of $V$ (in $B$)~; it is obtained by adjoining
$\root p\of u$ to $V$ for every $u\in\ogoth_V^\times$, where
$\ogoth_V$ is the ring of integers of $V$.  As $B=B'(\root p\of\pi)$
for every uniformiser $\pi$ of $K$, the group $\Gal(B|B')$ is cyclic
of order~$p$, the short exact sequence
$$
\{1\}\to\Gal(B|B')\to\Gal(B|K)\to\Gal(B'|K)\to\{1\}
$$
of profinite groups splits, and the resulting conjugation action of
$\Gal(B'|K)$ on $\Gal(B|B')$ is given by the cyclotomic character
$\omega:\Gal(B'|K)\to\F_p^\times$.  It follows
from \citer\iwasawa(Lemma 4) that the short exact sequence
$$
\{1\}\to\Gal(B'|V)\to\Gal(B'|K)\to\Gamma\to\{1\}
$$
of profinite groups also splits.  As the $\F_p[[\Gamma]]$-module
 $\Gal(B'|V)$ is isomorphic to
 $\Hom(\ogoth_V^\times\!/\ogoth_V^{\times p},\,{}_p\!V^\times)$, it is
 sufficient to determine the structure of the $\F_p[[\Gamma]]$-module
 $\ogoth_V^\times\!/\ogoth_V^{\times p}$, which we do.

\numeroter  Similarly, if $K$ has characteristic~$p$,
 then the short exact sequence
$$
\{1\}\to\Gal(B|V)\to\Gal(B|K)\to\Gamma\to\{1\}
$$
of profinite groups splits.  As the $\F_p[[\Gamma]]$-module
$\Gal(B|V)$ is isomorphic to $\Hom(V^+\!/\wp(V^+),\F_p)$, it is
sufficient to determine the structure of $V^+\!/\wp(V^+)$, which is
done below.

\numeroter Let $V_0$ be the maximal unramified extensions of $K$
(in $V$).  For every $n>0$, put $e_n=q^n-1$, $K_n=K(\!\root e_n\of1)$ and
$L_n=K_n(\!\root e_n\of{K_n^\times})$.  Note that $L_n$ is the maximal abelian
 extension of $K_n$ of exponent dividing~$e_n$, so it is galoisian
 over $K$~; put $G_n=\Gal(L_n|K)$.  The ramification index (resp.~the
 residual degree) of $L_n$ over $K$ is $e_n$ (resp.~$ne_n$).  We have
$$
V_0=\lim\limits_{\longrightarrow}K_n,\quad
V=\lim\limits_{\longrightarrow}L_n,\quad
\Gamma=\lim\limits_{\longleftarrow} G_n.
$$
Note that if $K$ has characteristic~$0$, then the $p$-torsion subgroup
${}_pL_n^\times$ of $L_n^\times$ has order~$p$ (because $L_1$ contains
$\root{p-1}\of{-p}\,$).

\numeroter Assume that $K$ has characteristic~$0$.  For
every finite extension $L$ of $K$, denote by $e_L$ ramification index
of $L|\Q_p$.  As $e_{L_n}\equiv0\pmod{e_n}$, we have
$e_{L_n}\equiv0\pmod{p-1}$, for every $n>0$~; write
$e_{L_n}=c_n.(p-1)$.  We have seen ($\oldstyle20$) that the
$\F_p[G_n]$-module $M_n=\bar U_{L_n}^1/\bar U_{L_n}^{c_n.p}$ is
isomorphic to $k[G_n]^{e_K}$.  Also, for every multiple $m$ of $n$,
the map $M_n\to M_m$ induced by the inclusion $L_n\subset L_m$ is {\it
injective}.  As $\ogoth_V^\times\!/\ogoth_V^{\times
p}=\lim\limits_{\longrightarrow}M_n$, we get from ($\oldstyle21$) by
passage to the limit~:

\numeroter {\it Let\/ $K$ be a finite extension of\/ $\Q_p$ of residue
field\/ $k$ and ramification index\/ $e_K$, let\/ $V$ be the maximal tamely
ramified extension of\/ $K$, and let\/ $\Gamma=\Gal(V|K)$.  The\/
$\F_p[[\Gamma]]$-module\/ $\ogoth_V^\times\!/\ogoth_V^{\times p}$ is
isomorphic to\/ $k[[\Gamma]]^{e_K}$, and the\/
$\F_p[[\Gamma]]$-module\/ $V^\times\!/V^{\times p}$ is isomorphic to
$k[[\Gamma]]^{e_K}\oplus\F_p$.}  \cqfd

\numeroter As for the dual $\F_p[[\Gamma]]$-modules
$\Gal(B'|V)=\Hom(\ogoth_V^\times\!/\ogoth_V^{\times p},{}_pV^\times)$
and $\Gal(B|V)=\Hom(V^\times\!/V^{\times p},{}_pV^\times)$, they are
respectivley isomorphic to $k[[\Gamma]]^{e_K}$ and to ${}_pV^\times
\oplus k[[\Gamma]]^{e_K}$.  Note that $k[[\Gamma]]^{e_K}$ is
free of rank $[K:\Q_p]$ over $\F_p[[\Gamma]]$.

\numeroter Now suppose that $K$ has characteristic~$p$.  By
an entirely similar argument, working with the modules
$M_n=\overline{L_n^+}/\overline{\ogoth_{L_n}^+}$
(resp.~$\overline{\ogoth_{L_n}^\times}$,
resp.~$\overline{L_n^\times}$), one gets from ($\oldstyle29$) and
($\oldstyle32$) by passage to the limit~:

\numeroter {\it For\/ $K=k\series{t}$, the\/ $\F_p[[\Gamma]]$-module\/ 
$\overline{V^+}=V^+/\wp(V^+)$ is isomorphic to\/
$\F_p[[\Gamma]]^{(\N)}$, and the\/ $\F_p[[\Gamma]]$-modules\/
$\overline{\ogoth_V^\times}=\ogoth_V^\times\!/\ogoth_V^{\times p}$
and\/ $\overline{V^\times}=V^\times\!/V^{\times p}$ are isomorphic
to\/ $\F_p[[\Gamma]]^{\N}$.}  \cqfd

\numeroter  As a result, the $\F_p[[\Gamma]]$-module
$\Gal(B|V)=\Hom(V^+\!/\wp(V^+),\F_p)$ is isomorphic to
$\F_p[[\Gamma]]^{\N}$.

\bigbreak
{\bf 6. Coronidis loco}
\bigskip

\numeroter  Let $K$ be a $p$-field and let $\tilde K$ be a 
maximal galoisian extension of $K$.  It is clear that if $K$ has
characteristic~$p$, then the profinite group $\Gal(\tilde K|K)$ cannot
be finitely generated, because $K$ has infinitely many cyclic
extensions of degree~$p$~: the dimension of the $\F_p$-space
$K^+\!/\wp(K^+)$ is infinite.  It is common knowledge that if $K$ is a
finite extension of $\Q_p$, then $\Gal(\tilde K|K)$ can be generated
by $[K:\Q_p]+3$ elements, cf.~\citer\jannsen(p.~65).  As a small gift
for the reader who has made it so far, we indicate how the foregoing
can be used to give a nice little proof~; it relies on the following
observation about profinite groups.

\numeroter  We say that a subset $S$ of a profinite
group $G$ {\it generates\/} $G$ if $G$ is the only {\it closed\/}
subgroup of $G$ containing $S$.  A {\it finite\/} subset $\Pi$ of a
pro-$p$-group $P$ generates $P$ if and only if its image $\bar\Pi$ in
the maximal commutative quotient $\bar P$ of $P$ of exponent
dividing~$p$ generates $\bar P$ (Burnside's ``\thinspace
basis\thinspace'' theorem)~; see \citer\kochpext(4.10).

\numeroter Consider a short exact sequence $\{1\}\to P\to G\to\Delta\to\{1\}$ of
profinite groups such that $P$ is a pro-$p$-group (and a closed
subgroup of $G$), so that $\bar P$ is an $\F_p[[\Delta]]$-module.
Presumably, {\it if\/ $\Pi\subset P$ is a finite subset whose image
in\/ $\bar P$ generates the\/ $\F_p[[\Delta]]$-module\/ $\bar P$, and
if\/ $D\subset G$ is a finite subset whose image in $\Delta$
generates\/ $\Delta$, then their union $\Pi\cup D$ generates\/ $G$}.
This should follow from an argument similar to the one
in \citer\jannsen(Lemma~3.3).  In our application below, the extension
$G$ of $\Delta$ by $P$ splits.


\numeroter  If this presumption is true, then 
($\oldstyle39$) provides a simple proof of the fact that the profinite
group $G=\Gal(\tilde K|K)$ is generated by $[K:\Q_p]+3$ elements.
Indeed, take $P=\Gal(\tilde K|V)$, so that $G/P=\Gamma$ and $\bar
P=\Gal(B|V)$.  We know from Iwasawa \citer\iwasawa(Theorem~2) that
$\Gamma$ is generated by {\it two\/} elements (and moreover the
extension $G$ of $\Gamma$ by $P$ splits).  We have seen
($\oldstyle39$) that the $\F_p[[\Gamma]]$-module $\Gal(B|V)$ is
generated by $[K:\Q_p]+1$ elements.  Hence, if ($\oldstyle45$) holds,
then $G$ is generated by $[K:\Q_p]+3$ elements.

\bigbreak
\unvbox\bibbox 

\bye